\documentclass{article}
\usepackage{amsfonts}

\usepackage{amsmath}

\newtheorem{theorem}{Theorem}

\newtheorem{lemma}[theorem]{Lemma}

\newenvironment{proof}[1][Proof]{\textbf{#1.} }{\ \rule{0.5em}{0.5em}}
\input{tcilatex}

\begin{document}

\title{CGMY and Meixner Subordinators are Absolutely Continuous with respect
to One Sided Stable Subordinators. }
\author{Dilip B. Madan \\
Robert H. Smith School of Business\\
Van Munching Hall\\
University of Maryland\\
College Park, MD 20742 \and Marc Yor \\
Laboratoire de probabilit\'{e}s et Modeles al\'{e}atoires\\
Universit\'{e} Pierre et Marie Curie\\
4, Place Jussieu F 75252 Paris Cedex}
\date{August 16 2005\\
}
\maketitle

\begin{abstract}
We describe the CGMY and Meixner processes as time changed Brownian motions.
The CGMY uses a time change absolutely continuous with respect to the
one-sided stable $(Y/2)$ subordinator while the Meixner time change is
absolutely continuous with respect to the one sided stable $(1/2)$
subordinator$.$ The required time changes may be generated by simulating the
requisite one-sided stable subordinator and throwing away some of the jumps
as described in Rosinski (2001).
\end{abstract}

\section{Introduction}

L\'{e}vy processes are increasingly being used to model the local motion of
asset returns, permitting the use of distributions that are both skewed and
capable of matching the high levels of kurtosis observed in factors driving
equity returns. By way of examples we cite the normal inverse Gaussian
process (Barndorff-Nielsen (1998)), the hyperbolic process (Eberlein, Keller
and Prause (1998)), and the variance gamma process (Madan, Carr and Chang
(1998)). For the valuation of structured equity products the importance of
skewness is well recognized and has led to the development of local L\'{e}vy
processes (See Carr, Geman, Madan and Yor (2004)) that preserve skews in
forward implied volatility curves. It is also understood from the steepness
of implied volatility curves that tail events have significantly higher
prices than those implied by a Gaussian distribution with the consequence
that pricing distributions display high levels of excess kurtosis.

On a single asset one may simulate the L\'{e}vy process calibrated to the
prices of vanilla options to value equity structured products written on a
single underlier. Such a simulation (See Rosinski (2001)) may approximate
the small jumps using a diffusion process with the large jumps simulated as
a compound Poisson process where one uses the normalized large jump L\'{e}vy
measure as the density of jump magnitudes with the integral of the L\'{e}vy
measure over the large jumps serving as the jump arrival rate. However,
increasingly one sees multiasset structures being traded and this requires a
modeling of asset correlations. Given marginal L\'{e}vy processes one could
accomodate correlations if one can represent the L\'{e}vy process as time
changed Brownian motion. In this case we correlate the simulated processes
by correlating the Brownian motions while preserving the independent time
changes for each of the marginal underliers.

It is therefore useful to have representations of L\'{e}vy processes as time
changed Brownian motions. For some L\'{e}vy processes, like the variance
gamma process or the normal inverse Gaussian process, these are known by
construction of the L\'{e}vy process via such a representation. For other L%
\'{e}vy processes, like the $CGMY$ process (Carr, Geman, Madan and Yor
(2002), see also Koponen (1995), Boyarchenko and Levendorskii (1999, 2000))
or the $Meixner$ process (Schoutens and Teugels (1998) see also Gregelionis
(1999), Schoutens (2000), and Pitman and Yor (2003)), the process is defined
directly by its L\'{e}vy measure and it is not clear a priori whether the
processes can be represented as time changed Brownian motions. With a view
to enhancing the applicability of these processes, particularly with respect
to multiasset structured products, we develop the representations of these
processes as time changed Brownian motions.

Section 2 presents for completeness, some preliminary results on L\'{e}vy
processes that we employ in the subsequent development. In section 3 we
develop the $CGMY$ process as a time changed Brownian motion with drift,
where the law of the time change is absolutely continuous over finite time
intervals with respect to that of the one sided stable $Y/2$ subordinator.
The simulation of $CGMY$ as time changed Brownian motion is described in
section 3. Section 4 develops the time change for the Meixner process as
absolutely continuous with respect to the one-sided stable $1/2$
subordinator. Simulation strategies for the Meixner process based on these
representations are described in Section 5. Section 6 reports on the
simulation results using chi-squared goodness of fit tests. Section 7
concludes.

\section{Preliminary results on L\'{e}vy processes}

We present three results from the theory of L\'{e}vy processes that we make
critical use of in our subsequent development. The first result relates the L%
\'{e}vy measure of a process obtained on subordinating a Brownian motion to
the L\'{e}vy measure of the subordinator. The second result establishes a
criterion for the absolute continuity of a subordinator with respect to
another subordinator. The third result presents the detailed relationship
between the standard presentation of the characteristic function of a two
sided jump and one-sided jump stable L\'{e}vy process and its L\'{e}vy
measure. \ These are presented in three short subsections.

\subsection{L\'{e}vy measure of a subordinated Brownian motion}

Suppose the L\'{e}vy process $X(t)$ is obtained by subordinating Brownian
motion with drift (i.e. the process $\theta u+W(u),$ for $(W(u),u\geq 0)$ a
Brownian motion) by an independent subordinator $Y(t)$ with L\'{e}vy measure 
$\nu (dy).$ Then applying Sato (1999) theorem 30.1 we get that the L\'{e}vy
measure of the process $X(t)$ is given by $\mu (dx)$ where%
\begin{equation}
\mu (dx)=dx\int_{0}^{\infty }\nu (dy)\frac{1}{\sqrt{2\pi y}}e^{-\frac{%
(x-\theta y)^{2}}{2y}}.  \label{Sato1}
\end{equation}

\subsection{Absolute Continuity Criterion for subordinators}

Suppose we have two subordinators $T_{A}=(T_{A}(t),t\geq
0),T_{B}=(T_{B}(t),t\geq 0).$ The law of the subordinator $T_{A}$ is
absolutely continuous with respect to the subordinator $T_{B},$ on finite
time intervals, just if there exists a function $f(t)$ such that the L\'{e}%
vy measures $\nu _{A}(dt),\nu _{B}(dt)$ for the processes $T_{A}$ and $T_{B}$
respectively are related by 
\begin{equation}
\nu _{A}(dt)=f(t)\nu _{B}(dt)  \label{ac1}
\end{equation}%
and furthermore, (Sato (1999) Theorem 33.1)%
\begin{equation}
\int_{0}^{\infty }\nu _{B}(dt)\left( \sqrt{f(t)}-1\right) ^{2}<\infty .
\label{ac2}
\end{equation}

\subsection{Stable Processes}

The Stable L\'{e}vy process $\mathcal{S}(\sigma ,\alpha ,\beta )=(X(t),t\geq
0)$ with parameters $(\sigma ,\alpha ,\beta )$ ( For details see DuMouchel
(1973, 1975), Bertoin (1996), Samorodnitsky and Taqqu (1998) Nolan (2001),
Ito (2004) ) has a characteristic function in standard form%
\begin{equation*}
E[e^{iuX(t)}]=\exp (-t\Psi (u))
\end{equation*}%
where the characteristic exponent $\Psi (u)$ is given by 
\begin{eqnarray}
\Psi (u) &=&\sigma ^{\alpha }|u|^{\alpha }\left( 1-i\beta sign(u)\tan \left( 
\frac{\pi \alpha }{2}\right) \right) ,\text{ }\alpha \neq 1  \label{stablecf}
\\
&=&\sigma |u|\left( 1+i\beta sign(u)\frac{2}{\pi }\log (|u|)\right) ,\text{ }%
\alpha =1.  \notag
\end{eqnarray}%
The parameters satisfy the restrictions, $\sigma >0,0<\alpha \leq 2$ and $%
-1\leq \beta \leq 1.$ The one sided jump stable processes result when $\beta
=1$ and there are only positive jumps or $\beta =-1$ in which case there are
only negative jumps.

The L\'{e}vy density of the stable process is of the form 
\begin{equation}
k(x)=\frac{c_{p}}{x^{1+\alpha }}\mathbf{1}_{x>0}+\frac{c_{n}}{|x|^{1+\alpha }%
}\mathbf{1}_{x<0}  \label{stablelm1}
\end{equation}%
and we have that 
\begin{equation}
\beta =\frac{c_{p}-c_{n}}{c_{p}+c_{n}}.  \label{stablebeta}
\end{equation}

It remains to express $\sigma $ in terms of the parameters of the L\'{e}vy
measure. In the one sided case with only positive jumps we have 
\begin{equation}
\sigma =\left[ \frac{c_{p}\Gamma \left( \frac{\alpha }{2}\right) \Gamma
\left( 1-\frac{\alpha }{2}\right) }{2\Gamma (1+\alpha )}\right] ^{\frac{1}{%
\alpha }}  \label{stablesigma1}
\end{equation}%
and more generally for the two sided jump case we have 
\begin{equation}
\sigma =\left[ \frac{c_{p}+c_{n}}{2}\frac{\Gamma \left( \frac{\alpha }{2}%
\right) \Gamma \left( 1-\frac{\alpha }{2}\right) }{\Gamma (1+\alpha )}\right]
^{\frac{1}{\alpha }}.  \label{stablesigma2}
\end{equation}

Conversely, $c_{p}$ and $c_{n}$ may be computed in terms of $\beta $ and $%
\sigma .$

\section{CGMY as time changed Brownian motion}

We wish to write the $CGMY$ process in the form 
\begin{equation*}
X_{CGMY}(t)=\theta Y(t)+W(Y(t))
\end{equation*}%
for an increasing time change process given by a subordinator $(Y(t),t\geq
0) $ independent of the Brownian motion $(W(s),s\geq 0)$ .

The characteristic function of the $CGMY$ process is 
\begin{equation*}
E\left[ \exp \left( iuX_{CGMY}(t)\right) \right] =\left( \phi
_{CGMY}(u)\right) ^{t}=\exp \left( tC\Gamma (-Y)\left[ 
\begin{array}{c}
\left( M-iu\right) ^{Y}-M^{Y}+ \\ 
\left( G+iu\right) ^{Y}-G^{Y}%
\end{array}%
\right] \right)
\end{equation*}%
The complex exponentiation is defined via the complex logarithm with a
branch cut on the negative real axis with polar coordinate arguments for the
complex logarithm restricted to the interval $]-\pi ,+\pi ]$. The $CGMY$
process is defined as a pure jump L\'{e}vy process by its L\'{e}vy measure%
\begin{equation*}
k_{CGMY}(x)=C\left[ \frac{\exp (-G|x|)}{|x|^{1+Y}}\mathbf{1}_{x<0}+\frac{%
\exp \left( -Mx\right) }{x^{1+Y}}\mathbf{1}_{x>0}\right] .
\end{equation*}

On the other hand we have, in all generality, by conditioning on the time
change that 
\begin{eqnarray*}
E\left[ e^{iu\left( \theta Y(t)+W(Y(t)\right) }\right] &=&E\left[ \exp
\left( iu\theta Y(t)-\frac{Y(t)}{2}u^{2}\right) \right] \\
&=&E\left[ \exp \left( -\left( \frac{u^{2}}{2}-iu\theta \right) Y(t)\right) %
\right]
\end{eqnarray*}

Take $u(\lambda )$ to be any solution of 
\begin{equation*}
\lambda =\left( \frac{u^{2}}{2}-iu\theta \right) ;
\end{equation*}%
Then we have the Laplace transform of the time change subordinator as%
\begin{equation*}
E[e^{-\lambda Y(t)}]=\exp \left( tC\Gamma (-Y)\left[ \left( M-iu(\lambda
)\right) ^{Y}-M^{Y}+\left( G+iu(\lambda )\right) ^{Y}-G^{Y}\right] \right)
\end{equation*}

The solutions for $u$ are: 
\begin{equation*}
u=i\theta \pm \sqrt{2\lambda -\theta ^{2}}
\end{equation*}%
where we suppose that $\theta ^{2}<2\lambda .$

We shall see that a good choice for $\theta $ , for sufficiently large $%
\lambda ,$ is 
\begin{equation*}
\theta =\frac{G-M}{2}
\end{equation*}%
and in this case%
\begin{eqnarray*}
M-iu &=&\frac{G+M}{2}+i\sqrt{2\lambda -\left( \frac{G-M}{2}\right) ^{2}} \\
G+iu &=&\frac{G+M}{2}-i\sqrt{2\lambda -\left( \frac{G-M}{2}\right) ^{2}}.
\end{eqnarray*}

It follows that the Laplace transform of the subordinator is

\begin{eqnarray*}
E[e^{-\lambda Y(t)}] &=&\exp \left( tC\Gamma (-Y)\left[ 2r^{Y}\cos (\eta
Y)-M^{Y}-G^{Y}\right] \right) \\
r &=&\sqrt{2\lambda +GM} \\
\eta &=&\arctan \left( \frac{\sqrt{2\lambda -\left( \frac{G-M}{2}\right) ^{2}%
}}{\left( \frac{G+M}{2}\right) }\right)
\end{eqnarray*}

In the special case of $G=M$ we have 
\begin{equation*}
E[e^{-\lambda Y(t)}]=\exp \left( 2tC\Gamma (-Y)\left[ \left( 2\lambda
+M^{2}\right) ^{Y/2}\cos \left( Y\arctan \left( \frac{\sqrt{2\lambda }}{M}%
\right) \right) -M^{Y}\right] \right)
\end{equation*}

\subsection{The explicit time change for CGMY}

We shall show that the time change subordinator $Y(t)$ associated with the $%
CGMY$ process is absolutely continuous with respect to the one-sided stable $%
Y/2$ subordinator \ and in particular that its L\'{e}vy measure $\nu (dy)$
takes the form%
\begin{eqnarray}
\nu (dy) &=&\frac{K}{y^{1+\frac{Y}{2}}}f(y)dy  \notag \\
f(y) &=&e^{-\frac{(B^{2}-A^{2})y}{2}}E\left[ e^{-\frac{B^{2}y}{2}\frac{%
\gamma _{Y/2}}{\gamma _{1/2}}}\right]  \label{CGMYTCLM} \\
B &=&\frac{G+M}{2}  \notag \\
K &=&\left[ \frac{C\Gamma \left( \frac{Y}{4}\right) \Gamma \left( 1-\frac{Y}{%
4}\right) }{2\Gamma (1+\frac{Y}{2})}\right]  \notag
\end{eqnarray}%
where $\gamma _{\frac{Y}{2}},\gamma _{\frac{1}{2}\text{ }}$ are two
independent gamma variates with unit scale parameters and shape parameters $%
Y/2,1/2$ respectively. Further we explicitly evaluate the expectation in
equation (\ref{CGMYTCLM}) in terms of the Hermite functions as follows.%
\begin{equation*}
E\left[ e^{-\frac{B^{2}y}{2}\frac{\gamma _{Y/2}}{\gamma _{1/2}}}\right] =%
\frac{\Gamma \left( \frac{Y}{2}+\frac{1}{2}\right) }{\Gamma (Y)\Gamma (\frac{%
1}{2})}2^{Y}\left( \frac{B^{2}y}{2}\right) ^{\frac{Y}{2}}I\left( Y,B^{2}y,%
\frac{B^{2}y}{2}\right)
\end{equation*}%
where%
\begin{equation*}
I(\nu ,a,\lambda )=\int_{0}^{\infty }x^{\nu -1}e^{-ax-\lambda
x^{2}}dx=(2\lambda )^{-\nu /2}\Gamma (\nu )h_{-\nu }\left( \frac{a}{\sqrt{%
2\lambda }}\right)
\end{equation*}%
and $h_{-\nu }(z)$ is the Hermite function with parameter $-\nu $ (see e.g
Lebedev (1972), p 290-291).

\subsection{Determining the time change for CGMY}

For an explicit evaluation of the time change we begin by writing the $CGMY$
L\'{e}vy density in the form 
\begin{equation*}
k_{CGMY}(x)=C\frac{e^{Ax-B|x|}}{x^{1+Y}},\text{ where: }A=\frac{G-M}{2};%
\text{ }B=\frac{G+M}{2}
\end{equation*}

Henceforth, when we encounter a L\'{e}vy measure $\mu (dx)$ that is
absolutely continuous with respect to Lebesgue measure we shall denote its
density by $\mu (x).$ We now employ the result (\ref{Sato1}) and seek to
find a L\'{e}vy measure of a subordinator satisfying 
\begin{eqnarray*}
C\frac{e^{Ax-B|x|}}{|x|^{1+Y}} &=&\int_{0}^{\infty }\nu (dy)\frac{1}{\sqrt{%
2\pi y}}e^{-\frac{(x-\theta y)^{2}}{2y}} \\
&=&\int_{0}^{\infty }\nu (dy)\frac{1}{\sqrt{2\pi y}}e^{-\frac{x^{2}}{2y}-%
\frac{\theta ^{2}y}{2}+\theta x}
\end{eqnarray*}

We set $\theta =A$ and observe that the right choice for $\theta $ is $%
(G-M)/2$ as remarked earlier, and identify $\nu (dy)$ such that 
\begin{equation}
C\frac{e^{-B|x|}}{|x|^{1+Y}}=\int_{0}^{\infty }\nu (dy)\frac{1}{\sqrt{2\pi y}%
}e^{-\frac{x^{2}}{2y}-\frac{\theta ^{2}y}{2}}  \label{eq1}
\end{equation}

We now recognize that the L\'{e}vy measure for the $CGMY$ is (taking $C=%
\frac{\Gamma \left( \frac{Y}{2}\right) \Gamma \left( 1-\frac{Y}{2}\right) }{%
\Gamma \left( 1+\frac{Y}{2}\right) },$ now), that of the symmetric stable $Y$
L\'{e}vy process with L\'{e}vy measure tilted as 
\begin{equation*}
k_{CGMY}(x)=e^{Ax-B|x|}k_{Stable(Y)}(x).
\end{equation*}%
\bigskip We also know that 
\begin{equation*}
X_{Stable(Y)}(t)=B_{Y^{0}(t)}
\end{equation*}%
where $Y^{0}(t)$ is the one sided stable $Y/2$ subordinator, independent of
the Brownian motion $(B_{u})$ .

We now write 
\begin{equation*}
X_{CGMY}(t)=\theta Y^{(1)}(t)+W_{Y^{(1)}(t)}
\end{equation*}%
and we seek to relate the L\'{e}vy measures $\nu ^{(1)}$ and $\nu ^{(0)}$ of
the processes $Y^{(1)}$ and $Y^{(0)}.$

From the result (\ref{Sato1}) we may write 
\begin{eqnarray*}
\mu _{0}(x) &=&\int_{0}^{\infty }\nu ^{(0)}(dy)\frac{e^{-\frac{x^{2}}{2y}}}{%
\sqrt{2\pi y}} \\
\mu _{1}(x) &=&\int_{0}^{\infty }\nu ^{(1)}(dy)\frac{e^{-\frac{(x-\theta
y)^{2}}{2y}}}{\sqrt{2\pi y}}
\end{eqnarray*}%
Hence we must have that%
\begin{equation*}
\int_{0}^{\infty }\nu ^{(1)}(dy)\frac{e^{-\frac{(x-\theta y)^{2}}{2y}}}{%
\sqrt{y}}=e^{Ax-B|x|}\int_{0}^{\infty }\nu ^{(0)}(dy)\frac{e^{-\frac{x^{2}}{%
2y}}}{\sqrt{y}}
\end{equation*}%
Taking $\theta =A,$ we get:%
\begin{equation*}
\int_{0}^{\infty }\nu ^{(1)}(dy)\frac{e^{-\frac{x^{2}}{2y}-\frac{A^{2}y}{2}}%
}{\sqrt{y}}=e^{-B|x|}\int_{0}^{\infty }\nu ^{(0)}(dy)\frac{e^{-\frac{x^{2}}{%
2y}}}{\sqrt{y}}
\end{equation*}

We now use the well known fact that 
\begin{equation*}
e^{-B|x|}=\int_{0}^{\infty }du\frac{B}{\sqrt{2\pi u^{3}}}e^{-\frac{B^{2}}{2u}%
-\frac{x^{2}}{2}u}
\end{equation*}%
to write 
\begin{equation*}
\int_{0}^{\infty }\nu ^{(1)}(dy)\frac{e^{-\frac{x^{2}}{2y}-\frac{A^{2}y}{2}}%
}{\sqrt{y}}=\int_{0}^{\infty }du\frac{B}{\sqrt{2\pi u^{3}}}e^{-\frac{B^{2}}{%
2u}}\int_{0}^{\infty }\nu ^{(0)}(dy)\frac{e^{-\frac{x^{2}}{2}\left( \frac{1}{%
y}+u\right) }}{\sqrt{y}}
\end{equation*}

By uniqueness of Laplace transforms we get that for every function $f:%
\mathbb{R}^{+}\rightarrow \mathbb{R}^{+}$%
\begin{equation*}
\int_{0}^{\infty }\nu ^{(1)}(dy)\frac{e^{-\frac{A^{2}y}{2}}}{\sqrt{y}}%
f\left( \frac{1}{y}\right) =\int_{0}^{\infty }du\frac{B}{\sqrt{2\pi u^{3}}}%
e^{-\frac{B^{2}}{2u}}\int_{0}^{\infty }\nu ^{(0)}(dy)\frac{1}{\sqrt{y}}%
f\left( \frac{1}{y}+u\right)
\end{equation*}%
or equivalently that, for every function $g:\mathbb{R}_{+}\rightarrow 
\mathbb{R}_{+}$ 
\begin{eqnarray*}
\int_{0}^{\infty }\nu ^{(1)}(dy)\frac{e^{-\frac{A^{2}y}{2}}}{\sqrt{y}}g(y)
&=&\int_{0}^{\infty }du\frac{B}{\sqrt{2\pi u^{3}}}e^{-\frac{B^{2}}{2u}%
}\int_{0}^{\infty }\nu ^{(0)}(dy)\frac{1}{\sqrt{y}}g\left( \frac{y}{1+uy}%
\right) \\
&=&\int_{0}^{\infty }du\frac{B}{\sqrt{2\pi u^{3}}}e^{-\frac{B^{2}}{2u}%
}\int_{0}^{\frac{1}{u}}d\left( \frac{s}{1-us}\right) \frac{\nu ^{(0)}(\frac{s%
}{1-us})}{\sqrt{\frac{s}{1-us}}}g(s) \\
&=&\int_{0}^{\infty }du\frac{B}{\sqrt{2\pi u^{3}}}e^{-\frac{B^{2}}{2u}%
}\int_{0}^{\frac{1}{u}}\frac{ds}{(1-su)^{2}}\frac{\nu ^{(0)}(\frac{s}{1-us})%
}{\sqrt{\frac{s}{1-us}}}g(s)
\end{eqnarray*}

Hence it is the case that 
\begin{eqnarray*}
\nu ^{(1)}(y)e^{-\frac{A^{2}y}{2}} &=&\int_{0}^{\frac{1}{y}}\frac{duBe^{-%
\frac{B^{2}}{2u}}\nu ^{(0)}(\frac{y}{1-uy})}{\sqrt{2\pi (u(1-uy))^{3}}} \\
&=&\sqrt{y}\int_{0}^{1}\frac{dvBe^{-\frac{B^{2}y}{2v}}\nu ^{(0)}(\frac{y}{1-v%
})}{\sqrt{2\pi (v(1-v))^{3}}}
\end{eqnarray*}%
In particular we have 
\begin{equation*}
\nu ^{(1)}(y)=\sqrt{y}\int_{0}^{1}\frac{dvBe^{-\frac{y}{2}\left( \frac{B^{2}%
}{v}-A^{2}\right) }\nu ^{(0)}(\frac{y}{1-v})}{\sqrt{2\pi (v(1-v))^{3}}}
\end{equation*}

We now introduce the explicit form of $\nu _{0}(y)$ for our case where it is
the L\'{e}vy density of the one-sided stable $Y/2$ subordinator, 
\begin{equation*}
\nu _{0}(y)=\frac{K}{y^{\left( \frac{Y}{2}+1\right) }}.
\end{equation*}%
This gives the representation%
\begin{eqnarray*}
\nu _{1}(y) &=&\frac{K}{y^{\frac{Y+1}{2}}}\int_{0}^{1}\frac{dvBe^{-\frac{y}{2%
}\left( \frac{B^{2}}{v}-A^{2}\right) }(1-v)^{\left( \frac{Y}{2}+1\right) }}{%
\sqrt{2\pi (v(1-v))^{3}}} \\
&=&\frac{K}{y^{\frac{Y+1}{2}}}\int_{1}^{\infty }\frac{dw}{w^{2}}\frac{Be^{-%
\frac{y}{2}\left( B^{2}w-A^{2}\right) }}{\sqrt{2\pi (\frac{1}{w}(1-\frac{1}{w%
}))^{3}}}\left( 1-\frac{1}{w}\right) ^{\left( \frac{Y}{2}+1\right) } \\
&=&\frac{K}{y^{\frac{Y+1}{2}}}\int_{1}^{\infty }\frac{dw}{\sqrt{2\pi w}}Be^{-%
\frac{y}{2}\left( B^{2}w-A^{2}\right) }\left( \frac{w-1}{w}\right) ^{\frac{%
Y-1}{2}} \\
&=&\frac{KBe^{-\frac{y}{2}\left( B^{2}-A^{2}\right) }}{y^{\frac{Y+1}{2}}}%
\int_{0}^{\infty }\frac{dh}{\sqrt{2\pi }}e^{-\frac{yB^{2}h}{2}}\frac{h^{%
\frac{Y-1}{2}}}{(1+h)^{\frac{Y}{2}}}
\end{eqnarray*}

\subsubsection{Absolute Continuity relations}

This subsection investigates the absolute continuity relation in general
between two subordinated processes and the absolute continuity of the
subordinators as processes. It is easy to show that the laws of the $CGMY$
process and the symmetric stable $Y$ process are locally equivalent, i.e.
for each $t,$their laws, as restricted to their past $\sigma -fields$ $%
\mathcal{F}_{t}$ up to time $t,$ are equivalent (from now on, as a slight
abuse of language, we shall say of 2 such processes, that they are
equivalent). Now that we have identified these processes as subordinated
processes, we look for the equivalence in law of the subordinators. Indeed
we first observe that if the subordinators are equivalent then the
subordinated processes will be equivalent but the converse may not be true.

Indeed, consider two subordinators 
\begin{equation*}
T_{A}(t),\text{ }T_{B}(t)
\end{equation*}%
such that the relation (\ref{ac1}) between their L\'{e}vy measures holds for
some function $f(t)$ for $t>0.$

We suppose the absolute continuity of $T_{A}$ with respect to $T_{B}$ or the
condition (\ref{ac2}).

We also define the subordinated processes 
\begin{eqnarray*}
X_{A}(t) &=&\beta _{T_{A}(t)} \\
X_{B}(t) &=&\beta _{T_{B}(t)}
\end{eqnarray*}%
where $(\beta _{u})$ is a Brownian motion assumed to be independent of
either $T_{A}$ or $T_{B}.$

We have from the result (\ref{Sato1}) that at the level of L\'{e}vy measures 
$\mu _{A},\mu _{B}$ for $X_{A},X_{B}$ 
\begin{eqnarray*}
\mu _{A}(x) &=&\int_{0}^{\infty }\nu _{A}(dt)\frac{e^{-\frac{x^{2}}{2t}}}{%
\sqrt{2\pi t}} \\
\mu _{B}(x) &=&\int_{0}^{\infty }\nu _{B}(dt)\frac{e^{-\frac{x^{2}}{2t}}}{%
\sqrt{2\pi t}}
\end{eqnarray*}

The following then holds as a consequence of (\ref{ac2}), for every
functional $F\geq 0$: 
\begin{equation*}
E\left[ F\left( T_{A}(s),s\leq t\right) \right] =E\left[ F(T_{B}(s),s\leq
t)\phi \left( T_{B}(s),s\leq t\right) \right]
\end{equation*}%
where 
\begin{equation*}
\phi \left( T_{B}(s),s\leq t\right) =\left( \frac{dP_{T_{A}}}{dP_{T_{B}}}%
\right) _{t}
\end{equation*}

As a consequence we deduce that, for every $G\geq 0:$ 
\begin{equation*}
E\left[ G\left( X_{A}(s),s\leq t\right) \right] =E\left[ G\left(
X_{B}(s),s\leq t\right) \phi (T_{B}(s),s\leq t)\right]
\end{equation*}

Consequently we may write%
\begin{equation*}
E\left[ G\left( X_{A}(s),s\leq t\right) \right] =E\left[ G\left(
X_{B}(s),s\leq t\right) \psi (X_{B}(s),s\leq t)\right]
\end{equation*}%
where%
\begin{equation*}
\psi (X_{B}(s),s\leq t)=E\left[ \phi (T_{B}(s),s\leq t)|(X_{B}(s),s\leq t)%
\right]
\end{equation*}

This implies that we should have 
\begin{equation*}
\mu _{A}(dx)=g(x)\mu _{B}(dx)
\end{equation*}%
with 
\begin{equation}
\int_{-\infty }^{\infty }\left( \sqrt{g(x)}-1\right) ^{2}\mu _{B}(dx)<\infty
\label{integ}
\end{equation}

We want to show that (\ref{ac2}) implies (\ref{integ}).

Now we have explicitly that 
\begin{eqnarray*}
g(x) &=&\frac{\int \nu _{A}(dt)\frac{e^{-\frac{x^{2}}{2t}}}{\sqrt{t}}}{\int
\nu _{B}(dt)\frac{e^{-\frac{x^{2}}{2t}}}{\sqrt{t}}} \\
&=&\frac{\int \nu _{B}(dt)f(t)\frac{e^{-\frac{x^{2}}{2t}}}{\sqrt{t}}}{\int
\nu _{B}(dt)\frac{e^{-\frac{x^{2}}{2t}}}{\sqrt{t}}}
\end{eqnarray*}

Let 
\begin{equation*}
\gamma ^{(x)}(dt)=\frac{\nu _{B}(dt)\frac{e^{-\frac{x^{2}}{2t}}}{\sqrt{t}}}{%
\int \nu _{B}(dt)\frac{e^{-\frac{x^{2}}{2t}}}{\sqrt{t}}}
\end{equation*}%
and note that 
\begin{equation*}
g(x)=\int \gamma ^{(x)}(dt)f(t)
\end{equation*}

We then have 
\begin{equation*}
\sqrt{g(x)}-1=\left( \int \gamma ^{(x)}(dt)f(t)\right) ^{\frac{1}{2}}-1
\end{equation*}

and 
\begin{equation*}
\int (\sqrt{g(x)}-1)^{2}\mu _{B}(dx)=\int (\sqrt{g(x)}-1)^{2}\left( \int \nu
_{B}(dt)\frac{e^{-\frac{x^{2}}{2t}}}{\sqrt{2\pi t}}dx\right)
\end{equation*}

Observe that 
\begin{eqnarray*}
&&(\sqrt{g(x)}-1)^{2}\mu _{B}(x) \\
&=&\left( \left( \int \gamma ^{(x)}(dt)f(t)\right) ^{\frac{1}{2}}-1\right)
^{2}\int \nu _{B}(dt)\frac{e^{-\frac{x^{2}}{2t}}}{\sqrt{2\pi t}} \\
&=&\left( \left( \frac{\int \nu _{B}(dt)f(t)\frac{e^{-\frac{x^{2}}{2t}}}{%
\sqrt{2\pi t}}}{\int \nu _{B}(dt)\frac{e^{-\frac{x^{2}}{2t}}}{\sqrt{2\pi t}}}%
\right) ^{\frac{1}{2}}-1\right) ^{2}\int \nu _{B}(dt)\frac{e^{-\frac{x^{2}}{%
2t}}}{\sqrt{2\pi t}} \\
&=&\int \nu _{B}(dt)f(t)\frac{e^{-\frac{x^{2}}{2t}}}{\sqrt{2\pi t}}+\int \nu
_{B}(dt)\frac{e^{-\frac{x^{2}}{2t}}}{\sqrt{2\pi t}} \\
&&-2\left( \int \nu _{B}(dt)f(t)\frac{e^{-\frac{x^{2}}{2t}}}{\sqrt{2\pi t}}%
\right) ^{\frac{1}{2}}\left( \int \nu _{B}(dt)\frac{e^{-\frac{x^{2}}{2t}}}{%
\sqrt{2\pi t}}\right) ^{\frac{1}{2}}
\end{eqnarray*}

We wish to show that the integral over $x$ of the right hand side is smaller
than 
\begin{equation*}
\int \nu _{B}(dt)f(t)+\int \nu _{B}(dt)-2\int \nu _{B}(dt)\sqrt{f(t)}
\end{equation*}

and this follows provided 
\begin{equation*}
\int dx\left( \int \nu _{B}(dt)f(t)\frac{e^{-\frac{x^{2}}{2t}}}{\sqrt{2\pi t}%
}\right) ^{\frac{1}{2}}\left( \int \nu _{B}(dt)\frac{e^{-\frac{x^{2}}{2t}}}{%
\sqrt{2\pi t}}\right) ^{\frac{1}{2}}\geq \int \nu _{B}(dt)\sqrt{f(t)}
\end{equation*}

For this consider 
\begin{eqnarray*}
&&\int \nu _{B}(dt)\sqrt{f(t)}=\int \nu _{B}(dt)\sqrt{f(t)}\int dx\frac{e^{-%
\frac{x^{2}}{2t}}}{\sqrt{2\pi t}} \\
&=&\int dx\int \nu _{B}(dt)\sqrt{f(t)}\left( \frac{e^{-\frac{x^{2}}{2t}}}{%
\sqrt{2\pi t}}\right) ^{\frac{1}{2}}\left( \frac{e^{-\frac{x^{2}}{2t}}}{%
\sqrt{2\pi t}}\right) ^{\frac{1}{2}} \\
&\leq &\int dx\left( \int \nu _{B}(dt)f(t)\frac{e^{-\frac{x^{2}}{2t}}}{\sqrt{%
2\pi t}}\right) ^{\frac{1}{2}}\left( \int \nu _{B}(dt)\frac{e^{-\frac{x^{2}}{%
2t}}}{\sqrt{2\pi t}}\right) ^{\frac{1}{2}},\text{ }
\end{eqnarray*}%
where we have used Cauchy-Schwarz for fixed $x.$

Hence we have 
\begin{equation*}
\int (\sqrt{g(x)}-1)^{2}\mu _{B}(dx)\leq \int_{0}^{\infty }\left( \sqrt{f(t)}%
-1\right) ^{2}\nu _{B}(dt)
\end{equation*}

The result does not go in the other direction as we may take 
\begin{eqnarray*}
\nu _{A}(dt) &=&\varepsilon _{a}(dt) \\
\nu _{B}(dt) &=&\varepsilon _{b}(dt)
\end{eqnarray*}%
for $a\neq b.$ These are not equivalent subordinators but in this case 
\begin{eqnarray*}
\mu _{A}(x) &=&\frac{e^{-\frac{x^{2}}{2a}}}{\sqrt{2\pi a}} \\
\mu _{B}(x) &=&\frac{e^{-\frac{x^{2}}{2b}}}{\sqrt{2\pi b}}
\end{eqnarray*}%
two L\'{e}vy densities, which in fact are probability densities, so that the
corresponding L\'{e}vy processes which are indeed Compound Poisson, are
(locally) equivalent.

\subsubsection{Absolute Continuity of the subordinators for CGMY and Stable
Y/2}

We now establish precisely the absolute continuity relationship between the
subordinator associated with the CGMY process, and the one sided stable $Y/2$
subordinator.

We note that 
\begin{eqnarray*}
\nu _{CGMY}(dy) &=&f(y)\nu _{0}(dy) \\
f(y) &=&e^{-\frac{y}{2}\left( B^{2}-A^{2}\right) }\left( B\sqrt{y}\right)
\int_{0}^{\infty }\frac{dh}{\sqrt{2\pi }}e^{-\frac{B^{2}y}{2}h}\frac{h^{%
\frac{Y-1}{2}}}{\left( 1+h\right) ^{\frac{Y}{2}}}
\end{eqnarray*}

We first check that as $B\rightarrow 0$ for $A=0$ we get the expected result
that $f(y)\rightarrow 1.$

For this we let $z=B\sqrt{y}$ and make the change of variable 
\begin{equation*}
k=z^{2}h
\end{equation*}%
to get 
\begin{eqnarray*}
f(y) &=&e^{-\frac{z^{2}}{2}}\int_{0}^{\infty }\frac{dk}{\sqrt{2\pi }z}e^{-%
\frac{k}{2}}\frac{\left( \frac{k}{z^{2}}\right) ^{\frac{Y-1}{2}}}{\left( 1+%
\frac{k}{z^{2}}\right) ^{\frac{Y}{2}}} \\
&=&e^{-\frac{z2}{2}}\int_{0}^{\infty }\frac{dk}{\sqrt{2\pi k}\frac{z}{\sqrt{k%
}}}e^{-\frac{k}{2}}\frac{\left( \frac{k}{z^{2}}\right) ^{\frac{Y-1}{2}}}{%
\left( 1+\frac{k}{z^{2}}\right) ^{\frac{Y}{2}}} \\
&\rightarrow &\int_{0}^{\infty }\frac{dk}{\sqrt{2\pi k}}e^{-\frac{k}{2}},%
\text{ as }z\rightarrow 0 \\
&=&2\int_{0}^{\infty }\frac{dx}{\sqrt{2\pi }}e^{-\frac{x^{2}}{2}} \\
&=&1
\end{eqnarray*}

For the equivalence of the two subordinators we must check that 
\begin{equation*}
\int_{0}^{\infty }\frac{dy}{y^{\frac{Y}{2}+1}}\left( \sqrt{f(y)}-1\right)
^{2}<\infty .
\end{equation*}

We break up this quantity in $2$ parts dealing with the integral near $0$
and $\infty $ separately. First consider the integral over $[1,\infty ).$
Here we write%
\begin{equation*}
\int_{1}^{\infty }\nu _{0}(dy)f(y)=\int_{1}^{\infty }\frac{dy}{y^{\frac{Y}{2}%
+1}}e^{-\frac{y}{2}(B^{2}-A^{2})}\left( B\sqrt{y}\right) \int_{0}^{\infty }%
\frac{dh}{\sqrt{2\pi }}e^{-\frac{B^{2}y}{2}h}\frac{h^{\frac{Y-1}{2}}}{\left(
1+h\right) ^{\frac{Y+2}{2}}}
\end{equation*}%
and check that 
\begin{equation*}
\left( B\sqrt{y}\right) \int_{0}^{\infty }\frac{dh}{\sqrt{2\pi }}e^{-\frac{%
B^{2}y}{2}h}\frac{h^{\frac{Y-1}{2}}}{\left( 1+h\right) ^{\frac{Y+2}{2}}}
\end{equation*}%
is bounded in $y.$

Write again $B\sqrt{y}=z,$ make the change of variable $k=z^{2}h$ and
observe that 
\begin{eqnarray*}
\left( B\sqrt{y}\right) \int_{0}^{\infty }\frac{dh}{\sqrt{2\pi }}e^{-\frac{%
B^{2}y}{2}h}\frac{h^{\frac{Y-1}{2}}}{\left( 1+h\right) ^{\frac{Y+2}{2}}}
&=&\int_{0}^{\infty }\frac{dk}{\sqrt{2\pi k}}e^{-\frac{k}{2}}\frac{\left( 
\frac{k}{z^{2}}\right) ^{\frac{Y}{2}}}{\left( 1+\frac{k}{z^{2}}\right) ^{%
\frac{Y}{2}}} \\
&\leq &\int_{0}^{\infty }\frac{dk}{\sqrt{2\pi k}}e^{-\frac{k}{2}}<\infty
\end{eqnarray*}

We next consider the required integral near $0,$ or over the interval $%
[0,1]. $ We have an expression of the form 
\begin{eqnarray*}
&&\int_{0}^{1}\frac{dy}{y^{\frac{Y}{2}+1}}\left( e^{-yC}\sqrt{I(y)}-1\right)
^{2}\text{ \ \ \ ,}C=\frac{B^{2}-A^{2}}{2} \\
I(y) &=&\left( B\sqrt{y}\right) \int_{0}^{\infty }\frac{dh}{\sqrt{2\pi }}e^{-%
\frac{B^{2}y}{2}h}\frac{h^{\frac{Y-1}{2}}}{\left( 1+h\right) ^{\frac{Y}{2}}}
\end{eqnarray*}

We now isolate the exponential by writing

\begin{eqnarray*}
&&\int_{0}^{1}\frac{dy}{y^{\frac{Y}{2}+1}}\left( e^{-yC}\sqrt{I(y)}-1\right)
^{2}=\int_{0}^{1}\frac{dy}{y^{\frac{Y}{2}+1}}\left( \left( e^{-yC}-1\right) 
\sqrt{I(y)}+\sqrt{I(y)}-1\right) ^{2} \\
&\leq &2\left( \int_{0}^{1}\frac{dy}{y^{\frac{Y}{2}+1}}\left(
e^{-yC}-1\right) ^{2}+\int_{0}^{1}\frac{dy}{y^{\frac{Y}{2}+1}}\left( \sqrt{%
I(y)}-1\right) ^{2}\right)
\end{eqnarray*}%
The exponential term is of order $y$ near zero and hence this first integral
is finite. For the second one we write 
\begin{eqnarray*}
&&\int_{0}^{1}\frac{dy}{y^{\frac{Y}{2}+1}}\left( \sqrt{I(y)}-1\right) ^{2} \\
&=&\int_{0}^{1}\frac{dy}{y^{\frac{Y}{2}+1}}\left[ \frac{\left( \sqrt{I(y)}%
-1\right) \left( \sqrt{I(y)}+1\right) }{\left( \sqrt{I(y)}+1\right) }\right]
^{2} \\
&\leq &\int_{0}^{1}\frac{dy}{y^{\frac{Y}{2}+1}}\left( I(y)-1\right) ^{2}
\end{eqnarray*}

For the finiteness of this integral we analyse the behavior of $(1-I(y))$
near zero. For this we analyse $I(y)=J(yB^{2})$ where 
\begin{eqnarray*}
J(y) &=&\sqrt{y}\int_{0}^{\infty }\frac{dh}{\sqrt{2\pi }}e^{-\frac{yh}{2}}%
\frac{h^{\frac{Y-1}{2}}}{\left( 1+h\right) ^{\frac{Y}{2}}} \\
&=&\int_{0}^{\infty }\frac{dk}{\sqrt{2\pi k}}e^{-\frac{k}{2}}\Phi \left( 
\frac{k}{y}\right) \\
&=&2\int_{0}^{\infty }\frac{dx}{\sqrt{2\pi }}e^{-\frac{x^{2}}{2}}\Phi \left( 
\frac{x^{2}}{y}\right) \\
&\equiv &\int_{-\infty }^{\infty }\frac{dx}{\sqrt{2\pi }}e^{-\frac{x^{2}}{2}%
}\Phi \left( \frac{x^{2}}{y}\right)
\end{eqnarray*}%
where%
\begin{equation*}
\Phi \left( \xi \right) =\left( \frac{\xi }{1+\xi }\right) ^{\frac{Y}{2}}
\end{equation*}

\begin{lemma}
The function $\Phi (\xi )$ is the distribution function of a random variable 
$V$ that can also be realized as the ratio of two independent gamma
variates, specifically 
\begin{equation*}
V\overset{(d)}{=}\frac{\gamma _{\frac{Y}{2}}}{\gamma \frac{1}{2}}
\end{equation*}%
where $\gamma _{a}$ is the gamma variate of parameter $a.$ In particular $V$
has finite moments of all orders $m<1$ and 
\begin{equation*}
E[V^{m}]=\frac{\Gamma \left( \frac{Y}{2}+m\right) }{\Gamma \left( \frac{Y}{2}%
\right) }\Gamma \left( 1-m\right)
\end{equation*}
\end{lemma}

\begin{proof}
\bigskip We note that $\Phi $ is the distribution function of a random
variable $V$ where for a uniform variate $U$ we have 
\begin{eqnarray*}
P\left( V\leq \xi \right) &=&P\left( U\leq \left( \frac{\xi }{1+\xi }\right)
^{\frac{Y}{2}}\right) \\
&=&P\left( U^{\frac{2}{Y}}\leq \frac{\xi }{1+\xi }\right) \\
&=&P\left( (1+\xi )U^{\frac{2}{Y}}\leq \xi \right) \\
&=&P\left( U^{\frac{2}{Y}}\leq \xi \left( 1-U^{\frac{2}{Y}}\right) \right) \\
&=&P\left( \frac{U^{\frac{2}{Y}}}{1-U^{\frac{2}{Y}}}\leq \xi \right)
\end{eqnarray*}

so that $V$ is the random variable 
\begin{equation*}
V=\frac{U^{\frac{2}{Y}}}{1-U^{\frac{2}{Y}}}
\end{equation*}

From the Beta-Gamma algebra we deduce that $V$ is 
\begin{equation*}
V=\frac{\gamma _{\frac{Y}{2}}}{\gamma _{1}}
\end{equation*}

Consequently $V$ has finite moments for all powers below unity. In
particular for $m<1$%
\begin{equation*}
E[V^{m}]=\frac{\Gamma \left( \frac{Y}{2}+m\right) }{\Gamma \left( \frac{Y}{2}%
\right) }\Gamma \left( 1-m\right)
\end{equation*}
\end{proof}

As a consequence for $m=\frac{1}{2}$ we have that%
\begin{equation*}
E[\sqrt{V}]=\frac{\Gamma \left( \frac{Y+1}{2}\right) }{\Gamma \left( \frac{Y%
}{2}\right) }\sqrt{\pi }.
\end{equation*}

Furthermore we have that as 
\begin{eqnarray*}
1-J(y) &=&P\left( |G|\leq \sqrt{Vy}\right) \\
&\sim &\sqrt{\frac{2}{\pi }}\sqrt{y}E\left[ \sqrt{V}\right]
\end{eqnarray*}

So the order of convergence of $1-I(y)=1-J(yB^{2})$ is always $\alpha =\frac{%
1}{2}$ and so 
\begin{equation*}
\frac{Y}{2}<2\alpha \equiv 1
\end{equation*}%
for all $Y<2.$ The desired absolute continuity result is established.

We also observe that 
\begin{equation*}
I(y)=J(yB^{2})=P\left( |G|\geq B\sqrt{Vy}\right) =P\left( \frac{G^{2}}{B^{2}V%
}\geq y\right) \leq 1
\end{equation*}

\subsubsection{A Further analysis of $I(y)$}

We now write the L\'{e}vy measure of the $CGMY$ subordinator in the form 
\begin{equation*}
\frac{K}{y^{1+\frac{Y}{2}}}E[e^{-yZ}]
\end{equation*}%
for some random variable $Z.$

For a fixed constant $B$ the L\'{e}vy measure of our subordinator in the
symmetric case is 
\begin{eqnarray*}
R &=&\frac{KBe^{-B^{2}\frac{y}{2}}}{y^{\frac{Y+1}{2}}}\int_{0}^{\infty }%
\frac{dh}{\sqrt{2\pi }}e^{-\frac{yB^{2}h}{2}}\frac{h^{\frac{Y-1}{2}}}{\left(
1+h\right) ^{\frac{Y}{2}}} \\
&=&\frac{KBe^{-B^{2}\frac{y}{2}}}{y^{\frac{Y+1}{2}}}\int_{0}^{\infty }\frac{%
dh}{\sqrt{2\pi h}}e^{-\frac{yB^{2}h}{2}}P\left( V\leq h\right) \\
&=&\frac{KBe^{-B^{2}\frac{y}{2}}}{y^{\frac{Y+1}{2}}}\int_{0}^{\infty }\frac{%
dk}{\sqrt{2\pi }}e^{-\frac{yB^{2}k^{2}}{2}}P\left( V\leq k^{2}\right) \\
&=&\frac{Ke^{-B^{2}\frac{y}{2}}}{y^{\frac{Y}{2}+1}}2\int_{0}^{\infty }\frac{%
dz}{\sqrt{2\pi }}e^{-\frac{z^{2}}{2}}P\left( V\leq \frac{z^{2}}{B^{2}y}%
\right) \\
&=&\frac{Ke^{-B^{2}\frac{y}{2}}}{y^{\frac{Y}{2}+1}}P\left( \frac{G^{2}}{%
VB^{2}}\geq y\right)
\end{eqnarray*}

We also know that 
\begin{equation*}
V\overset{(d)}{=}\frac{\gamma _{\frac{Y}{2}}}{\gamma _{1}}
\end{equation*}%
with two independent gamma variables. Thus we may write%
\begin{eqnarray*}
R &=&\frac{Ke^{-B^{2}\frac{y}{2}}}{y^{\frac{Y}{2}+1}}P\left( \gamma _{1}\geq 
\frac{yB^{2}\gamma _{\frac{Y}{2}}}{G^{2}}\right) \\
&=&\frac{Ke^{-B^{2}\frac{y}{2}}}{y^{\frac{Y}{2}+1}}E\left[ \exp \left( -%
\frac{yB^{2}\gamma _{\frac{Y}{2}}}{G^{2}}\right) \right] \\
&=&\frac{K}{y^{\frac{Y}{2}+1}}E\left[ \exp \left( -yB^{2}\frac{\gamma _{%
\frac{Y}{2}+\frac{y}{2}G^{2}}}{G^{2}}\right) \right]
\end{eqnarray*}

But 
\begin{equation*}
\frac{1}{2}G^{2}\overset{(d)}{=}\gamma _{\frac{1}{2}}
\end{equation*}%
so that we get 
\begin{equation}
R=\frac{Ke^{-B^{2}\frac{y}{2}}}{y^{\frac{Y}{2}+1}}E\left[ \exp \left( -y%
\frac{B^{2}}{2}\frac{\gamma _{\frac{Y}{2}}}{\gamma _{\frac{1}{2}}}\right) %
\right]  \label{CGMYsub}
\end{equation}

We now have identified the two L\'{e}vy measures as 
\begin{equation*}
\nu _{0}(dy)=\frac{Kdy}{y^{\frac{Y}{2}+1}}
\end{equation*}%
and 
\begin{eqnarray*}
\nu _{1}(dy) &=&\nu _{0}(dy)e^{-\frac{B^{2}y}{2}}E\left[ \exp \left(
-yZ\right) \right] \\
Z &=&\frac{B^{2}}{2}\frac{\gamma _{Y/2}}{\gamma _{1/2}}.
\end{eqnarray*}

\subsubsection{Evaluating explicitly the LT of Z}

There is an additional randomness in the simulation if the expectation 
\begin{equation*}
E[e^{-y\frac{B^{2}}{2}\frac{\gamma _{Y/2}}{\gamma _{1/2}}}]
\end{equation*}%
is evaluated by simulation. It is helpful to explicitly evaluate this
function. We begin with 
\begin{equation*}
\phi _{a,b}(\lambda )=E\left[ \exp \left( -\lambda \frac{\gamma _{a}}{\gamma
_{b}}\right) \right]
\end{equation*}

Now we have that 
\begin{eqnarray*}
e^{-\lambda }\phi _{a,b}(\lambda ) &=&E\left[ \exp \left( -\frac{\lambda }{%
\beta (b,a)}\right) \right] \\
&=&\frac{1}{B(b,a)}\int_{0}^{1}(1-x)^{a-1}x^{b-1}e^{-\frac{\lambda }{x}}dx \\
&=&\frac{1}{B(b,a)}\int_{1}^{\infty }\frac{dy}{y^{2}}\left( \frac{1}{y}%
\right) ^{b-1}\left( 1-\frac{1}{y}\right) ^{a-1}e^{-\lambda y} \\
&=&\frac{1}{B(b,a)}\int_{1}^{\infty }\frac{dy}{y^{a+b}}(y-1)^{a-1}e^{-%
\lambda y}
\end{eqnarray*}%
Hence we have that%
\begin{equation*}
\phi _{a,b}(\lambda )=\frac{1}{B(a,b)}\int_{0}^{\infty }\frac{z^{a-1}}{%
(1+z)^{a+b}}e^{-\lambda z}dz
\end{equation*}

We are interested in the case $a=\frac{Y}{2},$ $b=\frac{1}{2}$ and so we
write%
\begin{equation*}
\phi _{\frac{Y}{2},\frac{1}{2}}(\lambda )=\frac{1}{B\left( \frac{Y}{2},\frac{%
1}{2}\right) }\int_{0}^{\infty }dx\text{ }x^{\frac{Y}{2}-1}(1+x)^{-\frac{Y}{2%
}-\frac{1}{2}}e^{-\lambda x}
\end{equation*}

From Gradshetyn and Ryzhik (1995) (3.38) (7) Page 319 we have%
\begin{eqnarray*}
\int_{0}^{\infty }dx\text{ }x^{\frac{Y}{2}-1}(1+x)^{-\frac{Y}{2}-\frac{1}{2}%
}e^{-\lambda x} &=&2^{\frac{Y}{2}}\Gamma \left( \frac{Y}{2}\right) e^{\frac{%
\lambda }{2}}D_{-Y}\left( \sqrt{2\lambda }\right) \\
&=&2^{\frac{Y}{2}}\Gamma \left( \frac{Y}{2}\right) h_{-Y}\left( \sqrt{%
2\lambda }\right)
\end{eqnarray*}%
where $h_{\nu }(x)$ is the Hermite function of index $\nu $.

We have related the Hermite functions to the functions 
\begin{equation*}
I(\nu ,a,\lambda )=\int_{0}^{\infty }x^{\nu -1}e^{-ax-\lambda
x^{2}}dx=(2\lambda )^{-\nu /2}\Gamma (\nu )h_{-\nu }\left( \frac{a}{\sqrt{%
2\lambda }}\right)
\end{equation*}%
in \ Carr, Geman, Madan and Yor (2005).

We may therefore write%
\begin{equation*}
h_{-Y}(\sqrt{2\lambda })=\left( 2\lambda \right) ^{\frac{Y}{2}}\frac{1}{%
\Gamma (Y)}I(Y,2\lambda ,\lambda )
\end{equation*}%
It follows that%
\begin{eqnarray*}
\int_{0}^{\infty }dx\text{ }x^{\frac{Y}{2}-1}(1+x)^{-\frac{Y}{2}-\frac{1}{2}%
}e^{-\lambda x} &=&2^{\frac{Y}{2}}\Gamma \left( \frac{Y}{2}\right) \frac{%
\left( 2\lambda \right) ^{\frac{Y}{2}}}{\Gamma (Y)}I(Y,2\lambda ,\lambda ) \\
&=&2^{Y}\lambda ^{\frac{Y}{2}}\frac{\Gamma \left( \frac{Y}{2}\right) }{%
\Gamma (Y)}I(Y,2\lambda ,\lambda )
\end{eqnarray*}

It follows that 
\begin{equation*}
\phi _{\frac{Y}{2},\frac{1}{2}}(\lambda )=2^{Y}\lambda ^{\frac{Y}{2}}\frac{%
\Gamma \left( \frac{Y}{2}+\frac{1}{2}\right) }{\Gamma (Y)\Gamma (\frac{1}{2})%
}I\left( Y,2\lambda ,\lambda \right)
\end{equation*}

We therefore evaluate 
\begin{equation}
E\left[ e^{-y\frac{B^{2}}{2}\frac{\gamma _{\frac{Y}{2}}}{\gamma _{\frac{1}{2}%
}}}\right] =\frac{\Gamma \left( \frac{Y}{2}+\frac{1}{2}\right) }{\Gamma
(Y)\Gamma (\frac{1}{2})}2^{Y}\left( \frac{B^{2}y}{2}\right) ^{\frac{Y}{2}%
}I\left( Y,B^{2}y,\frac{B^{2}y}{2}\right)  \label{EXP}
\end{equation}

Putting together the result of equation (\ref{EXP}) and equation (\ref%
{CGMYsub}) we get the results for the $CGMY$ subordinator (\ref{CGMYTCLM}).

\section{Simulating CGMY using Rosinski Rejection}

We suppose that we have two L\'{e}vy measures $Q(dx),Q_{0}(dx)$ with the
property that 
\begin{equation*}
\frac{dQ}{dQ_{0}}\leq 1;
\end{equation*}%
and this is our case, then it is shown in Rosinski that we may simulate the
paths of $Q$ from those of $Q_{0}$ by only accepting all jumps $x$ in the
paths of $Q_{0}$ for which 
\begin{equation*}
\frac{dQ}{dQ_{0}}(x)>w
\end{equation*}%
where $w$ is an independent draw from a uniform distribution.

For our case we have that 
\begin{equation*}
\frac{d\nu _{1}}{d\nu _{0}}=E\left[ e^{-yZ}\right] <1
\end{equation*}%
and so accept all jumps in the paths of $\nu _{0}$ for which 
\begin{equation*}
E\left[ e^{-yZ}\right] >w
\end{equation*}

The detailed algorithm is for parameters $C,G,M,Y$ to first define the time
step to be $C,$%
\begin{equation*}
t=C.
\end{equation*}
Then we let 
\begin{eqnarray*}
A &=&\frac{G-M}{2} \\
B &=&\frac{G+M}{2}
\end{eqnarray*}

We next simulate at time $t$ from the one-sided stable subordinator with L%
\'{e}vy measure 
\begin{equation*}
\frac{1}{y^{\frac{Y}{2}+1}}dy
\end{equation*}

For this we let $\varepsilon =.0001$ and truncate jumps below $\varepsilon $
replacing them by their expected value at a rate of 
\begin{eqnarray*}
d &=&\int_{0}^{\varepsilon }y\frac{1}{y^{\frac{Y}{2}+1}}dy \\
&=&\frac{\varepsilon ^{1-\frac{Y}{2}}}{1-\frac{Y}{2}}
\end{eqnarray*}

For the arrival rate of jumps we have an arrival rate $\lambda $ of 
\begin{eqnarray*}
\lambda &=&\int_{\varepsilon }^{\infty }\frac{1}{y^{\frac{Y}{2}+1}}dy \\
&=&\frac{2}{Y}\frac{1}{\varepsilon ^{\frac{Y}{2}}}
\end{eqnarray*}

The interval jump times are exponential and are simulated by 
\begin{equation*}
t_{i}=-\frac{1}{\lambda }\log \left( 1-u_{2i}\right)
\end{equation*}%
for an independent uniform sequence $u_{2i}.$ The actual jump times are 
\begin{equation*}
\Gamma _{j}=\sum_{i=1}^{j}t_{i}
\end{equation*}

For the jump magnitude we simulate from the normalized L\'{e}vy measure the
jump size $y_{j}$ given by 
\begin{equation*}
y_{j}=\frac{\varepsilon }{\left( 1-u_{1j}\right) ^{\frac{2}{Y}}}
\end{equation*}%
for an independent uniform sequence $u_{1j}.$

The process $S(t)$ for the stable subordinator is given by 
\begin{equation*}
S(t)=dt+\sum_{j=1}^{\infty }y_{j}\mathbf{1}_{\Gamma _{j}<t}
\end{equation*}%
We now get the $CGMY$ subordinator $H(t)$ by 
\begin{eqnarray*}
H(t) &=&dt+\sum_{j=1}^{\infty }y_{j}\mathbf{1}_{\Gamma _{j}<t}\mathbf{1}%
_{h(y)>u_{3}{}_{j}} \\
h(y) &=&e^{-\frac{B^{2}y}{2}}\frac{\Gamma \left( \frac{Y}{2}+\frac{1}{2}%
\right) }{\Gamma (Y)\Gamma (\frac{1}{2})}2^{Y}\left( \frac{B^{2}y}{2}\right)
^{\frac{Y}{2}}I\left( Y,B^{2}y,\frac{B^{2}y}{2}\right)
\end{eqnarray*}%
for an independent uniform sequence $u_{3j}$

Finally we simulate the $CGMY$ random variable by 
\begin{equation*}
X=AH(t)+\sqrt{H(t)}z
\end{equation*}%
for a draw $z$ of a standard normal random variable.

\section{The Meixner Process as a Time Changed Brownian Motion}

We consider the Meixner Process (Schoutens and Teugels (1998), Pitman and
Yor (2003)) as a time changed Brownian motion. The L\'{e}vy measure of the
Meixner process is 
\begin{equation*}
k(x)=\delta \frac{\exp \left( \frac{b}{a}x\right) }{x\sinh \left( \frac{\pi x%
}{a}\right) }
\end{equation*}

The characteristic function is given by%
\begin{eqnarray*}
\phi _{Meixner}(u) &=&E[e^{iuX_{1}}] \\
&=&\left( \frac{\cos (b/2)}{\cosh (au-ib)/2)}\right) ^{2\delta }
\end{eqnarray*}

To see this process as a time changed Brownian motion we wish to identify $%
l(u)$ the L\'{e}vy measure of a subordinator such that 
\begin{eqnarray*}
k(x) &=&\int_{-\infty }^{\infty }\frac{1}{\sqrt{2\pi y}}\exp \left( -\frac{%
\left( x-Ay\right) ^{2}}{2y}\right) l(y)dy \\
&=&e^{Ax}\int_{-\infty }^{\infty }\frac{1}{\sqrt{2\pi y}}\exp \left( -\frac{%
x^{2}}{2y}-\frac{A^{2}y}{2}\right) l(y)dy
\end{eqnarray*}

Hence we set 
\begin{equation*}
A=\frac{b}{a}
\end{equation*}%
and seek to write%
\begin{equation}
\delta \frac{1}{x\sinh \left( \frac{\pi x}{a}\right) }=\int_{0}^{\infty }%
\frac{1}{\sqrt{2\pi y}}\exp \left( -\frac{x^{2}}{2y}-\frac{A^{2}y}{2}\right)
l(y)dy  \label{mxnr}
\end{equation}

We transform the left hand side of (\ref{mxnr}) as follows.

We recall that 
\begin{equation*}
\frac{Cx}{\sinh (Cx)}=E\left[ \exp \left( -\frac{x^{2}}{2}T_{C}^{(3)}\right) %
\right]
\end{equation*}%
where $T_{C}^{(3)}$ $=\inf \left\{ t|R_{t}^{(3)}=C\right\} $ for $%
R_{t}^{(3)} $ the $BES(3)$ process.

Then we write

\begin{eqnarray*}
\delta \frac{1}{x\sinh \left( \frac{\pi x}{a}\right) } &=&\frac{\delta
\left( \frac{\pi x}{a}\right) }{\left( \frac{\pi x^{2}}{a}\right) \sinh
\left( \frac{\pi x}{a}\right) } \\
&=&\frac{\delta a}{\pi }\frac{1}{x^{2}}E\left[ \exp \left( -\frac{x^{2}}{2}%
T_{C}^{(3)}\right) \right] \\
&=&\frac{\delta a}{\pi }\frac{1}{x^{2}}E\left[ \exp \left( -\frac{x^{2}C^{2}%
}{2}T_{1}^{(3)}\right) \right]
\end{eqnarray*}%
with $C=\frac{\pi }{a}.$ Denote by $\theta (h)dh$ the law of $T_{1}^{(3)}$ .
We may then write%
\begin{eqnarray*}
\delta \frac{1}{x\sinh \left( \frac{\pi x}{a}\right) } &=&\frac{\delta a}{%
\pi }\int_{0}^{\infty }\frac{du}{2}\exp \left( -\frac{x^{2}u}{2}\right) E%
\left[ \exp \left( -\frac{x^{2}C^{2}}{2}T_{1}^{(3)}\right) \right] \\
&=&\frac{\delta a}{2\pi }\int_{0}^{\infty }duE\left[ \exp \left( -\frac{x^{2}%
}{2}\left( u+C^{2}T_{1}^{(3)}\right) \right) \right] \\
&=&\frac{\delta a}{2\pi }\int_{0}^{\infty }du\int_{0}^{\infty }\theta
(t)dt\exp \left( -\frac{x^{2}}{2}(u+C^{2}t)\right) \\
&=&\frac{\delta a}{2\pi }\int_{0}^{\infty }du\int_{u}^{\infty }\frac{dv}{%
C^{2}}\exp \left( -\frac{x^{2}v}{2}\right) \theta \left( \frac{v-u}{C^{2}}%
\right) \\
&=&\frac{\delta a}{2\pi }\int_{0}^{\infty }dv\exp \left( -\frac{x^{2}v}{2}%
\right) \int_{0}^{v}\frac{du}{C^{2}}\theta \left( \frac{v-u}{C^{2}}\right) \\
&=&\frac{\delta a}{2\pi }\int_{0}^{\infty }dv\exp \left( -\frac{x^{2}v}{2}%
\right) \int_{0}^{\frac{v}{C^{2}}}dh\theta (h) \\
&=&\int_{0}^{\infty }dv\exp \left( -\frac{x^{2}v}{2}\right) \widehat{\theta }%
(v)
\end{eqnarray*}%
where%
\begin{eqnarray*}
\widehat{\theta }(v) &=&\frac{\delta a}{2\pi }\int_{0}^{\frac{v}{C^{2}}%
}\theta (h)dh \\
&=&\frac{\delta a}{2\pi }P\left( T_{1}^{(3)}\leq \frac{v}{C^{2}}\right) \\
&=&\frac{\delta a}{2\pi }P\left( Max_{t\leq \frac{v}{C^{2}}}R_{t}^{(3)}\geq
1\right)
\end{eqnarray*}

We recall that 
\begin{equation*}
T_{1}^{(3)}\overset{(law)}{=}\frac{1}{\left( \max_{t\leq
1}R_{t}^{(3)}\right) ^{2}}
\end{equation*}

We now transform the right hand side of (\ref{mxnr}) to write%
\begin{equation*}
\int_{0}^{\infty }\frac{1}{\sqrt{2\pi y}}\exp \left( -\frac{x^{2}}{2y}-\frac{%
A^{2}y}{2}\right) l(y)dy=\int_{0}^{\infty }\frac{1}{\sqrt{2\pi v^{3}}}\exp
\left( -\frac{x^{2}v}{2}-\frac{A^{2}}{2v}\right) l\left( \frac{1}{v}\right)
dv
\end{equation*}%
From the uniqueness of Laplace transforms we deduce that 
\begin{equation*}
\widehat{\theta }(v)=\frac{1}{\sqrt{2\pi v^{3}}}\exp \left( \frac{A^{2}}{2v}%
\right) l\left( \frac{1}{v}\right)
\end{equation*}%
or 
\begin{eqnarray*}
l(u) &=&\sqrt{\frac{2\pi }{u^{3}}}\widehat{\theta }\left( \frac{1}{u}\right)
\exp \left( -\frac{A^{2}u}{2}\right) \\
&=&\sqrt{\frac{2\pi }{u^{3}}}\frac{\delta a}{2\pi }P\left( M_{1}^{(3)}\geq C%
\sqrt{u}\right) \exp \left( -\frac{A^{2}u}{2}\right) \\
&=&\frac{\delta a}{\sqrt{2\pi u^{3}}}P\left( M_{1}^{(3)}\geq C\sqrt{u}%
\right) \exp \left( -\frac{A^{2}u}{2}\right) \\
&=&\frac{\delta a}{\sqrt{2\pi u^{3}}}g(u)
\end{eqnarray*}%
where%
\begin{equation*}
g(u)=P\left( M_{1}^{(3)}\geq C\sqrt{u}\right) \exp \left( -\frac{A^{2}u}{2}%
\right)
\end{equation*}

For the absolute continuity of our subordinator with respect to the one
sided stable $\frac{1}{2}$ subordinator we require that 
\begin{equation*}
\int \frac{1}{\sqrt{u^{3}}}\left( \sqrt{g(u)}-1\right) ^{2}du<\infty .
\end{equation*}%
For this we observe that 
\begin{eqnarray*}
\left( \sqrt{g(u)}-1\right) ^{2} &\leq &|g(u)-1| \\
&=&1-g(u) \\
&=&1-P\left( M_{1}^{(3)}\geq C\sqrt{u}\right) \exp \left( -\frac{A^{2}u}{2}%
\right) \\
&=&1-\exp \left( -\frac{A^{2}u}{2}\right) +\exp \left( -\frac{A^{2}u}{2}%
\right) \left( 1-P\left( M_{1}^{(3)}\geq C\sqrt{u}\right) \right)
\end{eqnarray*}%
The first part is clearly integrable with respect to $\left( \frac{du}{%
u^{3/2}}\right) $and for the second we observe that as%
\begin{equation*}
\lambda ^{k}P\left( T\geq \lambda \right) \leq E\left[ T^{k}\right]
\end{equation*}%
that%
\begin{equation*}
P\left( \frac{1}{(M_{1}^{(3)})^{2}}\geq \frac{1}{C^{2}u}\right) =P\left(
T_{1}^{(3)}\geq \frac{1}{C^{2}u}\right) \leq Ku^{k},\text{ for all }k
\end{equation*}

For the simulation of Meixner as a time changed Brownian motion we would
wish to evaluate 
\begin{eqnarray*}
P\left( M_{1}^{(3)}\geq C\sqrt{u}\right) &=&P\left( \frac{1}{%
(M_{1}^{(3)})^{2}}\leq \frac{1}{C^{2}u}\right) \\
&=&P\left( T_{1}^{(3)}\leq \frac{1}{C^{2}u}\right) \\
&=&P(\pi ^{2}T_{1}^{(3)}\leq \frac{\pi ^{2}}{C^{2}u}) \\
&=&P(T_{\pi }^{(3)}\leq \frac{\pi ^{2}}{C^{2}u}) \\
&=&\sum_{-\infty }^{\infty }(-1)^{n}e^{-n^{2}\pi ^{2}/(2C^{2}u)}
\end{eqnarray*}%
For the last equality we refer to Pitman and Yor (2003).

\section{Simulation of the Meixner Process}

The simulation strategy is similar to that employed in section 3 for $CGMY,$
except that here we simulate first the jumps of the one sided stable $\frac{1%
}{2}$ with L\'{e}vy density 
\begin{equation*}
k(x)=\frac{\delta a}{\sqrt{2\pi x^{3}}},\text{ }x>0.
\end{equation*}

We approximate the small jumps of the subordinator using the drift 
\begin{equation*}
\zeta =\delta a\sqrt{\frac{2\varepsilon }{\pi }}
\end{equation*}

The arrival rate for the jumps above $\varepsilon $ is 
\begin{equation*}
\lambda =\delta a\sqrt{\frac{2}{\pi \varepsilon }}
\end{equation*}%
and the jump sizes for the one sided stable$\left( \frac{1}{2}\right) $ are 
\begin{equation*}
y_{j}=\frac{\varepsilon }{u_{j}^{2}}
\end{equation*}%
for an independent uniform sequence $u_{j}.$

We then evaluate the function $g(y)$ at the point $y_{j}$ and define the
time change variable 
\begin{equation*}
\tau =\varsigma +\sum_{j}y_{j}\mathbf{1}_{g(y_{j})>w_{j}}
\end{equation*}%
for yet another independent uniform sequence $w_{j}.$ We note that the
function $g(y)$ only use the parameters $a,b$ and is independent of the
parameter $d.$

The value of the Meixner random variable or equivalently the unit time level
of the process is then generated as 
\begin{equation*}
X=\frac{b}{a}\tau +\sqrt{\tau }z
\end{equation*}%
where $z$ is an independent standard normal variate.

\section{Results of Simulations}

For both the $CGMY$ and $Meixner$ processes we present in this section the
results of simulating the processes at \ typical parameter values obtained
on calibrating option prices on the S\&P 500 index. The parameter values for
the $CGMY$ \ are $C=1,$ $G=5,$ $M=10,$ and $Y=.5.$ The parameters for the $%
Meixner$ were $a=.25,$ $b=-1.5$ and $\delta =1.$

We present graphs (\ref{tcbmcgmy},\ref{tcbmmxnr}) for a weekly time step $%
h=.02$ of the simulated and actual densities as well as chi square tests of
the hypothesis that the sample was drawn from the respective densities. The
solid lines are the theoretical density while the data points are indicated
by dots. The sample sizes in both cases were $5000.$ The range for both the $%
CGMY$ and $Meixner$ returns was $25\%.$ In both cases we used $100$ cells
and employed those with more than five observations for the test. The $CGMY$
had \ a chisquare statistic of $42.0122$ with $56$ degrees of freedom and a $%
p-value$ of $.9172.$ For the $Meixner$ the test statistic was $78.70$ with $%
84$ degrees of freedom and a $p-value$ of $.6427.$

\FRAME{ftbpFU}{5.0194in}{4.0145in}{0pt}{\Qcb{CGMY simulation as time changed
Brownian Motion using shaved one sided stable $Y/2.$}}{\Qlb{tcbmcgmy}}{%
tcbmcgmy.eps}{\special{language "Scientific Word";type
"GRAPHIC";maintain-aspect-ratio TRUE;display "USEDEF";valid_file "F";width
5.0194in;height 4.0145in;depth 0pt;original-width 6.749in;original-height
5.188in;cropleft "0";croptop "1.0390";cropright "1";cropbottom "0";filename
'../../MATLAB701/work/tcbmcgmy.eps';file-properties "XNPEU";}}

\bigskip

\FRAME{ftbpFU}{5.0194in}{4.0145in}{0pt}{\Qcb{Meixner simulation as time
changed Brownian motion using shaved one sided stable $1/2.$}}{\Qlb{tcbmmxnr}%
}{tcbmmxnr.eps}{\special{language "Scientific Word";type
"GRAPHIC";maintain-aspect-ratio TRUE;display "USEDEF";valid_file "F";width
5.0194in;height 4.0145in;depth 0pt;original-width 6.6642in;original-height
5.188in;cropleft "0";croptop "1.0261";cropright "1";cropbottom "0";filename
'../../MATLAB701/work/tcbmmxnr.eps';file-properties "XNPEU";}}

\bigskip

\pagebreak


\begin{thebibliography}{99}
\bibitem{BN} Barndorff-Nielsen, O.E. (1998), ``Processes of Normal Inverse
Gaussian type,''\textit{Finance and Stochastics, 2, 41-68.}

\bibitem{JB} Bertoin, J. (1996), \textit{L\'{e}vy Processes}, Cambridge
University Press, Cambridge.

\bibitem{BL} Boyarchenko, S.I. and Levendorskii (1999), ``Generalizations of
the Black-Scholes equation for Truncated L\'{e}vy processes,'' Working paper.

\bibitem{BL2} Boyarchenko, S.I. and Levendorskii (2000), ``Option pricing
for Truncated L\'{e}vy processes,'' \textit{International Journal for Theory
and Applications in Finance}, 3, 549-552.

\bibitem{cgmy} Carr, P. Geman, H., Madan, D. and M. Yor (2002), ``The Fine
Structure of Asset Returns: An Empirical Investigation,'' \textit{The
Journal of Business}, 75, 305-332.

\bibitem{cgmy2} Carr, P., Geman, H., Madan, D. and M. Yor (2004), ``From
Local Volatility to Local L\'{e}vy Models,'' \textit{Quantitative Finance},
5, 581-588.

\bibitem{cgmy3} Carr, P., Geman, H., Madan, D. and M. Yor (2005), ``Pricing
Options on Realized Variation,'' forthcoming in \textit{Finance and
Stochastics}

\bibitem{dumouchela} DuMouchel, W.H. (1973), ``Stable distribution in
statistical inference 1: Symmetric stable distributions compared to other
symmetric long-tailed distributions,'' \textit{Journal of the American
Statistical Association}, 68, 469-477.

\bibitem{dumouchelb} DuMouchel, W.H. (1975), ``Stable distributions in
statistical inference 2: Information from stably distributed samples,'' 
\textit{Journal of the American Statistical Association}, 70, 386-393.

\bibitem{EKP} Eberlein, E., U. Keller, and K. Prause (1998), ``New Insights
into Smile, Mispricing and Value at Risk,''\textit{Journal of Business}, 71,
371-406.

\bibitem{GR} Gradshetyn, I.S. and I.M. Ryzhik (1995), \textit{Table of
Integrals, Series and Products}, Academic Press, New York.

\bibitem{Grig} Grigelionis, B. (1999), ``Processes of Meixner Type,'' 
\textit{Lithuanian Mathematics Journal} 39, 33-41.

\bibitem{Ito} Ito, K. (2004), \textit{Stochastic Processes}, Springer,
Berlin.

\bibitem{Koponen} Koponen, I. (1995), ``Analytic Approach to the problem of
convergence of truncated L\'{e}vy flights towards the Gaussian stochastic
process,'' \textit{Physical Review} E52, 1197-1199.

\bibitem{Lebedev} Lebedev, N.N. (1972), \textit{Special Functions and Their
Applications}, Dover, New York.

\bibitem{MCC} Madan, D., P. Carr, and E. Chang (1998), ``The Variance Gamma
Process and Option Pricing,''\textit{European Finance Review}, 2, 79-105.

\bibitem{Nolan} Nolan. J.P. (2001), ``Maximum Likelihood Estimation and
Diagnostics for Stable Distributions,'' In \textit{L\'{e}vy processes-Theory
and Application}. O.E. Barndorff-Nielsen, T. Mikosch and S.I. Resnick Eds.,
Birkhauser, Boston, 379-400.

\bibitem{PY} Pitman, J. and M. Yor (2003), ``Infinitely Divisible Laws
associated with Hyperbolic Functions,'' Canadian Journal of Mathematics, 55,
2, 292-330.

\bibitem{rosinski} Ros\'{\i}nski, J. (2001), ``Series representations of L%
\'{e}vy processes from the perspective of point processes,'' In \textit{L%
\'{e}vy processes-Theory and Application}. O.E. Barndorff-Nielsen, T.
Mikosch and S.I. Resnick Eds., Birkhauser, Boston, 401-415.

\bibitem{Sato} Sato K. (1999), \textit{L\'{e}vy Processes and Infinitely
Divisible Distributions}, Cambridge University Press, Cambridge.

\bibitem{ST} Samorodnitsky, G. and M.S. Taqqu (1994), \textit{Stable
Non-Gaussian Random Processes}, Chapman and Hall, New York.

\bibitem{Sch} Schoutens, W. (2000), \textit{Stochastic Processes and
Orthogonal Polynomials}, Lecture Notes in Statistics 146, Springer, Berlin.

\bibitem{SchT} Schoutens, W. and J.L.Teugels (1998), ``L\'{e}vy processes,
polynomials and martingales,'' \textit{Communications in Statistics:
Stochastic Models}, 14, 335-349.
\end{thebibliography}
\end{document}